\documentclass[a4paper, reqno, 12pt]{amsart}

\usepackage{geometry}             

\usepackage{float}
\usepackage{bm}
\usepackage{fullpage,xcolor}
\usepackage[mathscr]{euscript}
\usepackage[all]{xy}
\usepackage{epsfig}
\usepackage{amsfonts}
\usepackage{mathptmx}

\usepackage{graphicx}
\usepackage{amssymb}
\usepackage{amsmath}
\usepackage{amsthm}
\usepackage{mathrsfs}
\usepackage{epstopdf}
\usepackage{url}
\usepackage{dirtytalk}
\usepackage[msc-links,alphabetic]{amsrefs}
\usepackage{tikz}

\usepackage{color}
\makeatletter

\@addtoreset{equation}{section}
\makeatother 

\theoremstyle{plain}
\newtheorem{theorem}{Theorem}[section]

\newtheorem{proposition}[theorem]{Proposition}
\newtheorem{lemma}[theorem]{Lemma}

\theoremstyle{definition}
\newtheorem{definition}[theorem]{Definition}

\newtheorem{remark}[theorem]{Remark}

\usepackage{latexsym}
 
\title[On the classification of periodic weaves and \\ universal cover of links in thickened surfaces] {On the classification of periodic weaves and \\ universal cover of links in thickened surfaces}

\author[Sonia Mahmoudi]{Sonia Mahmoudi}
\address{Advanced Institute for Materials Research, Tohoku University, 2-1-1 Katahira, Aoba-ku, Sendai, 980-8577, Japan}
\email{sonia.mahmoudi@tohoku.ac.jp}
\date{March 19, 2024}
\thanks{The author would like to thank the editor and the reviewers for their valuable comments. This work was supported by a Research Fellowship from JST CREST Grant Number JPMJCR17J4 and Grant-in Aid for JSPS Fellows Number 22J13397.} 

\begin{document}

\begin{abstract} 
A \textit{periodic weave} is the lift of a particular link embedded in a thickened surface of genus $g \geq 1$ to the universal cover. Its components are infinite unknotted simple open curves that can be partitioned in at least two distinct \textit{sets of threads}. The classification of periodic weaves can be reduced to the one of their generating cells, namely their \textit{weaving motifs}. However, this classification cannot be achieved through the classical theory of links in thickened surfaces since periodicity in the universal cover is not encoded. In this paper, we first introduce the notion of hyperbolic periodic weaves, which generalizes our doubly periodic weaves embedded in the Euclidean thickened plane. Then, Tait’s First and Second Conjectures are extended to minimal reduced alternating weaving motifs and proved using a generalized Kauffman bracket polynomial defined for periodic weaving diagrams of $\mathbb{E}^2$ and generalized to $\mathbb{H}^2$. The first conjecture states that any minimal alternating reduced weaving motif has the minimum possible number of crossings, while the second one formulates that two such oriented weaving motifs have the same writhe.   
\end{abstract}

\maketitle


\section{Introduction}\label{sec:1}

\begin{figure}[ht]
\centerline{\includegraphics[width=5in]{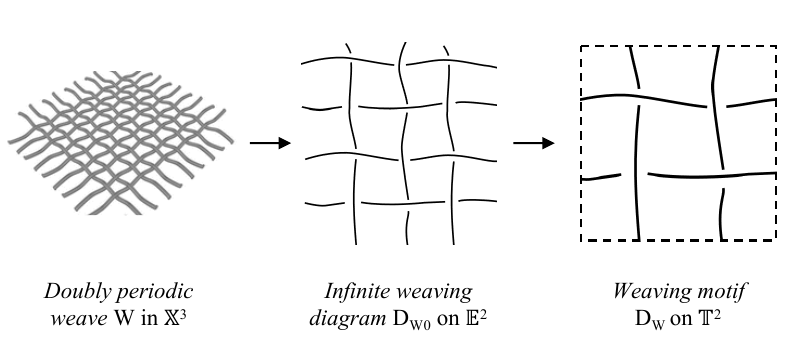}}
\vspace*{8pt}
\caption{\label{diagram} A doubly periodic weave (left), its infinite diagram (center), a weaving motif defined as the quotient by an integer lattice (right).}
\end{figure}

A \textit{periodic weave} is an entangled structure that can be described as the lift of a particular link embedded in a thickened surface of genus $g \geq 1$ to the Euclidean or hyperbolic thickened plane, denoted by $\mathbb{E}^2 \times I$ and $\mathbb{H}^2 \times I$ with $I=[-1,1]$, respectively. The topological study of doubly periodic weaves in $\mathbb{E}^2 \times I$ was first introduced by S. A. Grishanov et al. \cite{Grishanov1} in the context of classification of textile structures from a knot theoretical point of view. An example of doubly periodic weave is illustrated on the left of Figure~\ref{diagram}. Following this approach, various topological invariants have been constructed for the classification of doubly periodic textiles (\cite{GrishanovP1}, \cite{GrishanovP2}, \cite{Morton}, \cite{Kawauchi}, \cite{Vitally}). The common strategy used in these studies was to reduce the classification problem of doubly periodic structures to the one of their generating cells on the diagrammatic level, that we call \textit{motifs}. Note however that even if these motifs are defined as link diagrams on the torus $T^2$, their classification differs from the classical theory since the notion of equivalence classes is not defined similarly. For instance, a \textit{Dehn twist} of $T^2$, which corresponds to a change of basis in $\mathbb{E}^2$ equipped with a fixed point lattice or equivalently to a global isotopy (shearing) of the doubly periodic structure, induces an equivalence relation between motifs in addition to the usual isotopy moves for classical links \cite{Grishanov1}.

To distinguish weaves from other classes of entangled structures, we previously introduced a new definition of doubly periodic weaves embedded in $\mathbb{E}^2 \times I$. Our definition partitions the components of a weave, called \textit{threads}, in at least two disjoint \textit{sets of threads}, each being characterized by a `direction' (see \cite{Sonia1} and \cite{Sonia2}). For comparison, \textit{knits} (\cite{Kawauchi}) or \textit{braids} (\cite{Lambropoulou}) are defined by curve components running along a single direction. On the diagrammatic level, we call any generating cell of a doubly periodic weaving diagram a \textit{weaving motif} (see Figure~\ref{diagram} for an example).
In particular, a weaving motif is a link diagram whose components consist only of essential simple closed curves on the torus, which in other words are non-self-intersecting and non-contractible closed curves. Additionally, these curves lift to unknotted simple open curves in $\mathbb{E}^2 \times I$ and must belong to at least two distinct sets of threads, as recall in Section~\ref{sec:2-1}.

\begin{figure}[ht]
\centering
   \includegraphics[width=5in]{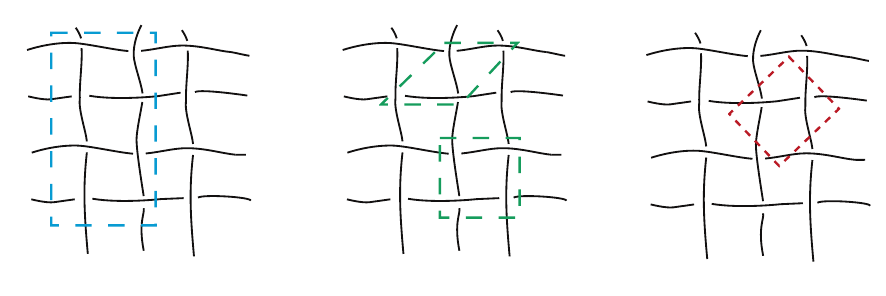}
      \caption{\label{crossing} Different generating cells for the same doubly periodic biaxial alternating weave. On the left, the number of crossings is 8 for the blue diagram.  In the center, two equivalent green diagrams with 4 crossings. On the right, a red minimal diagram with 2 crossings.}
\end{figure}

An important remark is now needed. By considering different point lattices, one may also obtain non-isotopic weaving motifs that represent the same periodic weave. For example in Figure~\ref{crossing}, the blue, green and red diagrams are weaving motifs which lift to the same periodic weave in $\mathbb{E}^2 \times I$. As mentioned above, note that for a given fixed integer lattice, one can select two generating cells randomly, simply by a change of basis as illustrated in the middle of Figure~\ref{crossing}. These two green motifs are not considered equivalent in the sense of the classical theory since they differ by a Dehn twist of $T^2$. As a notable consequence, note that two non-isotopic torus knots lift to the same periodic structure in $\mathbb{E}^2 \times I$. Therefore, we should consider that the classification of periodic weaves depends on the choice of the lattice. In this paper, we are in particular interested in classifying periodic weaves by invariants depending on the number of crossings. This leads to the need of introducing the notion of \textit{minimal lattice}, defined in Section~\ref{sec:2}, that generates what we should call \textit{minimal motifs}. Then, within an equivalence class of minimal motifs, each diagram that achieves the minimal number of crossings, namely the \textit{crossing number}, is said to be \textit{minimum}. Note that a minimum diagram is always minimal by definition but the converse may not hold. For example, in Figure~\ref{crossing}, the red motif is a minimal and minimum motif for the corresponding periodic weave. Besides, one can also consider that up to isotopy and torus twist, this red motif is contained in the green and blue motifs. We will thus introduce the notion of \textit{scale-equivalence}, to encode that the lifts to $\mathbb{E}^2 \times I$ of different finite covers of the same motif define the same doubly periodic weave in $\mathbb{E}^2 \times I$. These notions lead to the statement of our \textit{Generalized Reidemeister Theorem}~\ref{th:Doubly-R-Th} for doubly periodic weave in Section~\ref{sec:2-1}, which extends the one stated in \cite{Grishanov1}.

In this paper, we also generalize the notion of periodic weaves to the hyperbolic thickened plane $\mathbb{H}^2 \times I$ by introducing the notion of \textit{hyperbolic periodic weaves}. Hence, the corresponding infinite diagrams lie on $\mathbb{H}^2$, which we choose to represent by the Poincare disk model (see Section~\ref{sec:2-3}). 
The definition of weaving motifs as link diagrams on higher genus surfaces is obtained through the quotient of a periodic infinite diagram on $\mathbb{H}^2$ by a discrete group generated by translations. As for the case of doubly periodic weaves, the classical theory of links embedded in a thickened surface of genus $g >1$, which does not encode the periodicity in $\mathbb{H}^2 \times I$, fails to classify hyperbolic periodic weaves. In Section~\ref{sec:2-3}, we present a \textit{generalized Reidemeister theorem} for hyperbolic periodic weaves (Theorem~\ref{th:Periodic-R-Th}), which highlight the difference with the classical theory, since the notion of Dehn twists and scale-equivalence also extends for higher genus surfaces and their lift in $\mathbb{H}^2 \times I$.

\smallbreak

In the history of knot theory, much attention has been paid to classify \textit{alternating} links. An alternating link is defined as a link which admits at least one diagram whose crossings alternate between over and under as one travels around each component \cite{MurasugiBook}. In the late nineteenth century, P.G. Tait \cite{Tait} stated several famous conjectures on alternating links that remained unproven for a century, until the discovery of the Jones polynomial. The first conjecture classifies alternating links by their crossing number and was demonstrated independently for links in $\mathbb{S}^3$ by M.B Thistlethwaite \cite{Thistlethwaite}, K. Murasugi \cite{Murasugi1}, and L.H. Kauffman \cite{Kauffman}. More recently, Tait's first conjecture was generalized for classical links embedded in thickened surfaces by T. Fleming et al. \cite{AdamsFleming}, \cite{Fleming}, and H.U. Boden et al. \cite{Boden2022}, \cite{Boden2019}. Tait's second conjecture for alternating links follows from the first one and classifies oriented diagrams by their \textit{writhe}. The writhe of a link is the sum of the signs of all the crossings, where each crossing is assigned a sign $\pm1$. This second conjecture was originally proven independently by M.B Thistlethwaite \cite{Thistlethwaite} and K. Murasugi \cite{Murasugi2} for the case of links in $\mathbb{S}^3$. It has also been extended to classical links in a thickened surface by H.U. Boden et al. in \cite{Boden2022} and \cite{Boden2019}.

The notion of alternating links also appears at the level of weaves. In particular, an \textit{alternating} periodic weave can be defined as a weave that admits at least one alternating weaving motif. One can thus generalize Tait's conjectures to periodic weaves in $\mathbb{E}^2 \times I$ and $\mathbb{H}^2 \times I$. The main purpose of this paper is to prove these conjectures for alternating weaving motifs on a thickened surface of genus $g \geq 1$. As mentioned above, since the notion of equivalence differs from the classical theory, these two conjectures are proved with respect to our Generalized Reidemeister Theorems (Theorems~\ref{th:Doubly-R-Th},~\ref{th:Periodic-R-Th}), which encode the periodicity of weaves. 
More precisely, we will first prove the following.

\begin{theorem}\label{th:1-1} \textbf{(Tait’s First Conjecture for Alternating Weaving Motifs)}
A minimal reduced alternating weaving motif is a minimum diagram of its alternating periodic weave. 
\end{theorem}

To prove Theorem~\ref{th:1-1} we follow the combinatorial approach of \cite{Kauffman} presented for classic links in $\mathbb{S}^3$, where all steps are adapted for alternating periodic weaves. Note that the original proof indeed fails for the case of periodic weaves in $\mathbb{E}^2 \times I$ and $\mathbb{H}^2 \times I$ and the aim of this paper is to fill that gap. We do that by first defining the notion of \textit{reduced} weaving motif, which plays the role of a reduced diagram in the original proof. Intuitively, a reduced motif can be thought of as a diagram where no crossing can be removed, as detailed in Definition.~\ref{def:2-8}. Then, we use the generalized version of the \textit{bracket polynomial of weaving motifs} embedded in $T^2$ defined in \cite{GrishanovP2}, which applies for doubly periodic weaves in $\mathbb{E}^2 \times I$.  Besides, we also extend this bracket polynomial for weaving motifs on surfaces of genus $g >1$ so that Theorem~\ref{th:1-1} also applies for alternating hyperbolic periodic weaves. We do so by capturing the essential closed curve components of a state diagram into a $2g$-uplet of integers that is invariant under Dehn twists of the surface, as detailed in Section~\ref{sec:3}.

\smallbreak

Finally, we also naturally consider the generalization of Tait's second conjecture for periodic weaves. In particular, we will prove the following.

\begin{theorem}\label{th:1-2} \textbf{(Tait’s Second Conjecture for weaves)}
Any two connected minimal reduced alternating diagrams of an oriented doubly periodic weave have the same writhe.
\end{theorem}

The proof of Theorem.~\ref{th:1-2} follows the strategy of the one from R. Strong \cite{Stong} for classic links embedded in $\mathbb{S}^3$. 
The main difference when generalizing this results to periodic weaves is also the consideration of the choice of periodic lattices. In addition, one must also pay attention to the fact that the quotient by a point lattice of a component of a periodic weave, namely a thread, which is an unknotted open curve in the plane, can be a knotted closed curve on the corresponding motif. 

\bigbreak

This paper is organized as follows. In Section~\ref{sec:2}, we introduce definitions of periodic weaves and their diagrams, as well as the notions of alternating and reduced motifs. 
In Section~\ref{sec:3}, we first recall necessary results on the bracket polynomial for weaving motifs on a torus and we generalize this polynomial to motifs on higher genus surfaces using Teichm\"uller theory. In Section~\ref{sec:4}, we present the proof of the main theorem of this paper, that is Tait’s first conjecture for weaves. Finally, in Section~\ref{sec:5}, we expose the proof of Tait’s second conjecture for weaves, with the use of the proof of the first conjecture.

\newpage

\section{Periodic weaves and their corresponding links in thickened surfaces}\label{sec:2}

In this section, we first recall the definitions of doubly periodic weaves embedded in the Euclidean thickened plane and their corresponding weaving motifs in the thickened torus. Then, we generalize these notions to define periodic weaves embedded in the hyperbolic thickened plane and their corresponding weaving motifs in higher genus surfaces. Finally, we extend the notions of alternating and reduced link diagrams to the case of weaving motifs.

\subsection{Doubly periodic weaves and links in a thickened torus}\label{sec:2-1}

~

~

In this subsection, we recall the definition of a doubly periodic weave embedded in the Euclidean thickened plane as the lift of a particular type of link embedded in the thickened torus, that we presented in \cite{Sonia2}. 
\smallbreak

Let $\mathbb{E}^2$ denote the Euclidean plane, $T^2$ denote the torus and $I=[-1,1]$ be an interval. Let also $u,v$ be a basis of $\mathbb{E}^2$ such that the covering map $\rho: \mathbb{E}^2 \rightarrow{} \mathrm{T}^2$ sends $u$ and $v$ to the longitude $l$ and the meridian $m$ of $\mathrm{T}^2$, respectively. Note that this covering map extends to the thickened plane by $\tilde{\rho}: \mathbb{E}^2 \times I \to T^2 \times I$. Considering these notations, the set of points $\Lambda (u,v) = \{xu + yv / x,y \in \mathbb{Z}\}$ generated by the basis $u,v$ of $\mathbb{E}^2$ defines a point lattice of $\mathbb{E}^2$ isomorphic to $\mathbb{Z}^2$. 

Let now $t$ be an essential closed curve component of a link diagram on $T^2$, such that $t$ is homotopic to a $(p,q)$-torus knot. Recall that by a $(p,q)$-torus knot is meant an essential simple closed curve that intersects the torus meridian $p$ times and the torus longitude $q$ times, as detailed in \cite{MurasugiBook}(Chapter 7). Moreover, a $(p',q')$-torus link is a link in the torus consisting of $g$ components, each being a $(p,q)$-torus knot, where $p' = g \times p$ and $q' = g \times q$. Then, the lift of $t$ under the covering map $\rho$ is a set of infinite open curves related by a planar translation on $\mathbb{E}^2$, called \textit{threads}. We say that two threads belong to the \textit{same set of threads} and are \textit{parallel} if their quotient by $\Lambda (u,v)$ is either the same essential closed curve $t$ in $T^2$, or a set of two essential closed curve $t$ and $t'$ homotopic to the same $(p,q)$-torus knot in $T^2$. In particular, if two components $t$ and $t'$ of a link diagram in $T^2$ are homotopic to the same $(p,q)$-torus knot, they are also said to be \textit{parallel} and to belong to the \textit{same set of threads}.

Next, if $t$ admits a self-crossing on $T^2$, then we are interested in distinguishing the cases when this crossing vanishes on $\mathbb{E}^2$. In particular, if a thread has no self-crossing on $\mathbb{E}^2$, we say that it is \textit{unknotted}. To do so, we first introduce the notion of `divided curves' as follows. Recall that an essential closed curve component $t$ of a link diagram on $T^2$ lifts to open curves in $\mathbb{E}^2$. Moreover, $t$ is also said to be \textit{unknotted} if it does not admit any self-crossing, or \textit{knotted} if it admits a finite number of self-crossings, each being a transverse double point. Since both unknotted and knotted essential closed curves on $T^2$ may lift to unknotted open curves on $\mathbb{E}^2$, these two cases are considered in our definition of weaves. However, we want to exclude the case where a knotted component on $T^2$ lifts to knotted threads on $\mathbb{E}^2$. In particular, if $t$ is a knotted closed curve on $T^2$ with a self-crossing at a point $p$. Then, $p$ can divide $t$ into two closed curves, possibly knotted too. We refer to these two closed curves in $T^2$ as $t_p$ and $t'_p$, called \textit{divided curves} of $t$ for $p$. In particular, it is known that if at least one of the divided curves $t_p$ or $t'_p$ is null-homotopic on $T^2$, then $t$ lifts to knotted curves on $\mathbb{E}^2$. For details, we refer to \cite{Tanio}. 

\begin{figure}[ht]
\centerline{\includegraphics[width=5in]{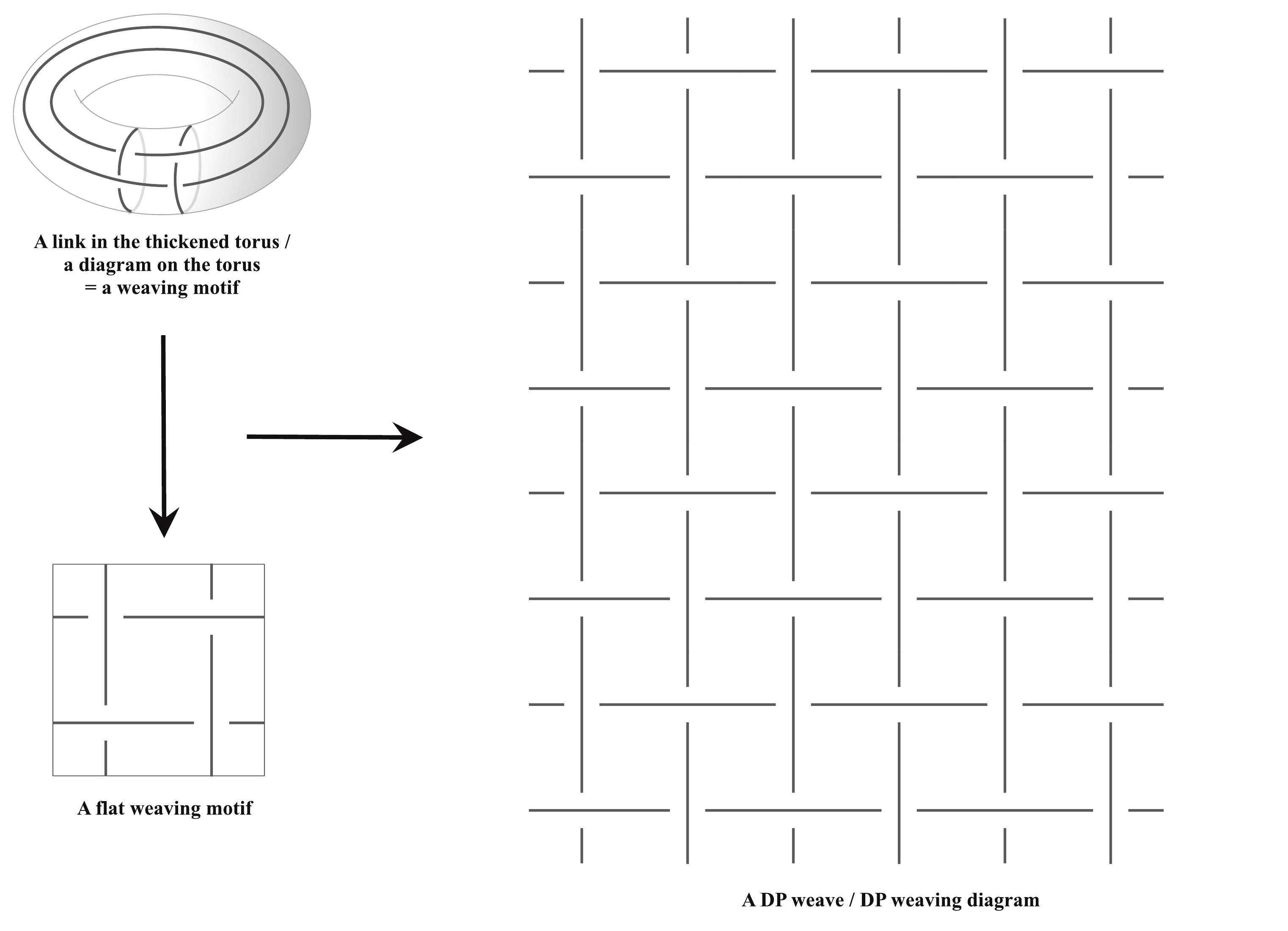}} 
\vspace*{8pt}
\caption{\label{DP-weave}
A weaving motif in the (thickened) torus (on the left), and its corresponding  doubly periodic weave (diagram) in the Euclidean (thickened) plane (on the right).}
\end{figure}

Doubly periodic weaves and their weaving motifs are thus defined as follows.
For an illustration, see Figure~\ref{DP-weave}.

\begin{definition}\label{def:EuclideanWeave} 
Let $W$ be a link embedded in $T^2 \times I$ with its corresponding diagram $D_W$ in $T^2$ satisfying the following conditions, with $i,j$ being positive integers, 
\begin{enumerate}
\item each component $t_i$ of $W$ is homotopic to a $(p_i,q_i)$-torus knot, where $p_i$ and $q_i$ are coprime integers,
 \smallbreak
\item each component $t_i$ of lifts to an unknotted thread under $\tilde{\rho}$,
 \smallbreak 
\item $W$ contains at least two distinct non-parallel components $t_i$ and $t_j$.
\end{enumerate}
\noindent Then, the lift of $W$ (resp. $D_W$) under $\tilde{\rho}$ (resp. $\rho$) to $\mathbb{E}^2 \times I$ (resp. $\mathbb{E}^2 $) is called a {\it doubly periodic weave}, denoted by $W_{\infty}$ (resp. a {\it doubly periodic weaving diagram}, denoted by $D_{\infty}$). Moreover, $D_W$ is said to be a {\it weaving motif} of $D_{\infty}$.
\end{definition}

Now recall that $\Lambda (u,v)$ defines a periodic lattice of unit parallelograms. 
In particular, each of these parallelograms can be identified with a flat torus on which is embedded a weaving motif, considering the identification space of $T^2$. 
In other words, a weaving motif can be defined as the quotient of a doubly periodic weaving diagram $D_{\infty}$ by $\Lambda (u,v)$.
However, it is well-known that different bases of $\mathbb{E}^2$ can generate equivalent point lattices. 
More specifically, 
$$\Lambda (u,v) \simeq \Lambda' (u',v') \, \textit{ if and only if } \, 
\begin{pmatrix}
u'  \\
v' 
\end{pmatrix}
=
\begin{pmatrix}
a & b \\
c & d 
\end{pmatrix}
\times
\begin{pmatrix}
u \\
v 
\end{pmatrix}
\textit{ , where }
\mid{ad - bc}\mid = \pm 1$$

\begin{figure}[ht]
\centerline{\includegraphics[width=5in]{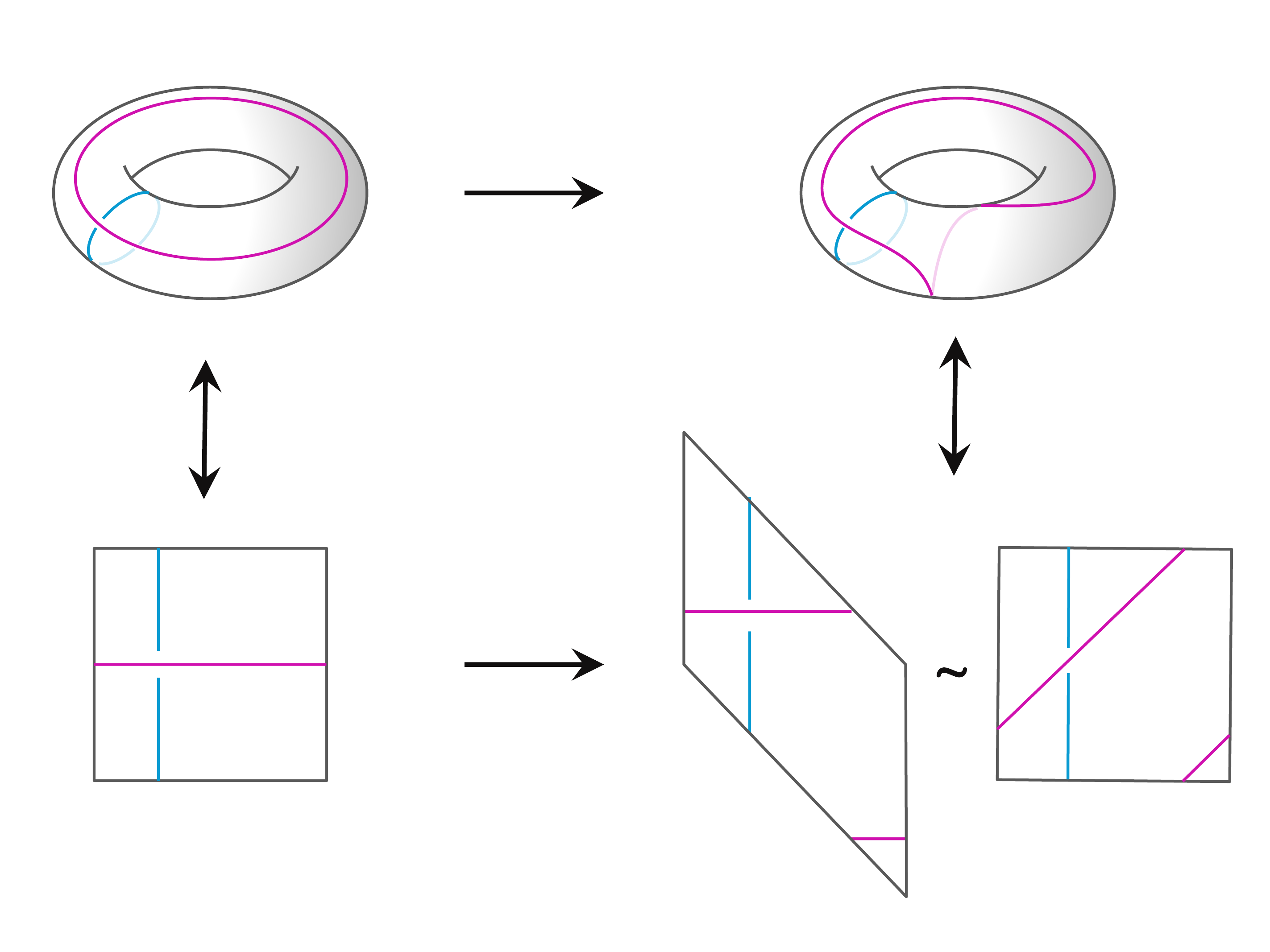}}
\vspace*{8pt}
\caption{\label{twist} 
A Dehn twist visualized on a flat torus.}
\end{figure}

This implies that different weaving motifs may lift to the same doubly periodic weaving diagram. It is well-known that this equivalence of point lattices can be translated at the level of the torus. Considering linear orientation-preserving homeomorphism of $\mathbb{E}^2$, namely matrices of $SL(2,\mathbb{Z})$, the equivalence of generating cells that differ by a change of basis of $\mathbb{E}^2$ can be translated into Dehn twists of $T^2$, see Figure~\ref{twist} for illustration and \cite{Farb}, \cite{Grishanov1} for details. It follows that for a fixed point lattice, the equivalence of doubly periodic weaves can be described combinatorially by the following.

\begin{proposition}\label{R-theorem} \textbf{(Generalized Reidemeister Theorem for Fixed Point Lattices in $\mathbb{E}^2$ \cite{Grishanov1})}
Two doubly periodic weaves in $\mathbb{E}^2 \times I$ are \textit{isotopic} for equivalent point lattices if and only if their corresponding weaving motifs can be obtained from each other by a finite sequence of \textit{Reidemeister moves} $R_1$, $R_2$, and $R_3$, torus isotopies and Dehn twists. 
\end{proposition}

\begin{figure}[ht]
\centerline{\includegraphics[width=5in]{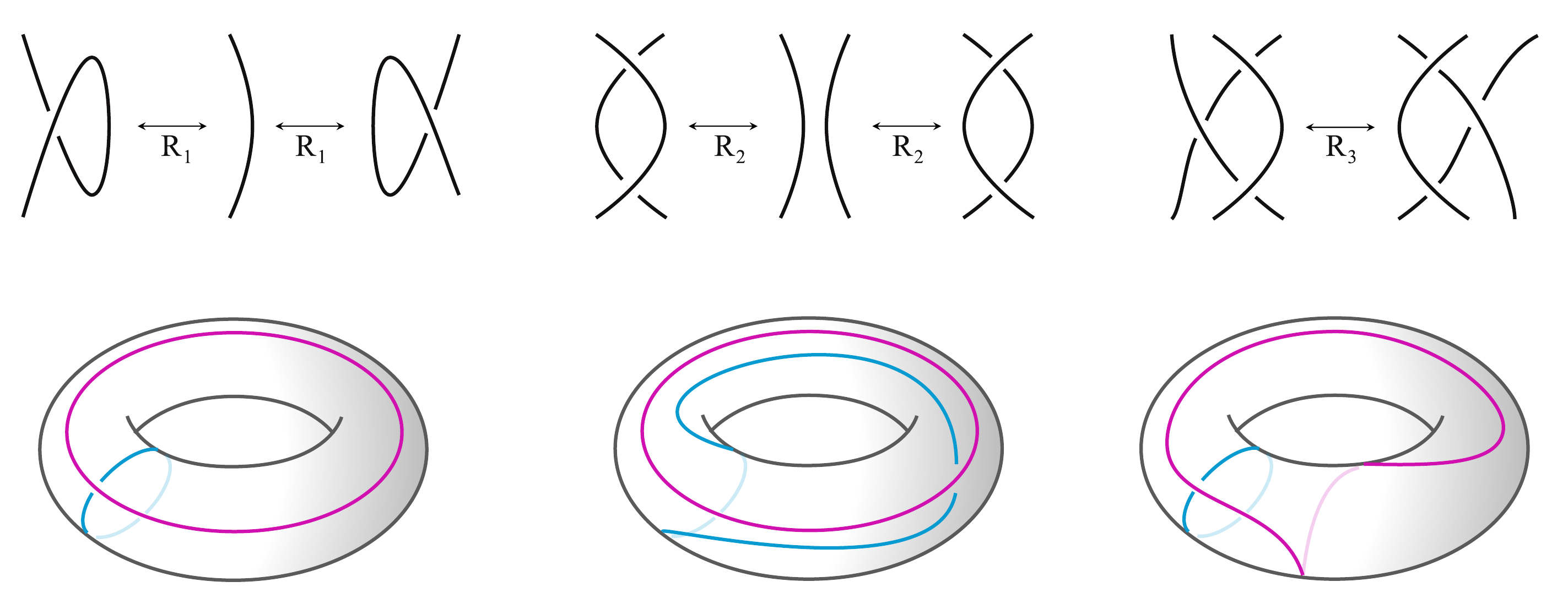}}
\vspace*{8pt}
\caption{\label{Rmoves} 
On the top, the three classical Reidemeister moves.
On the bottom, a torus with a meridian and a longitude curves (left), its longitudinal twist (middle) and its meridional twist (right)}
\end{figure}

However, this notion of equivalence is not strong enough to capture non-equivalent point lattices associated with weaving motifs that lift to equivalent periodic weaving diagrams, as illustrated in Figure~\ref{crossing}. Consider for example,  
$$\Lambda (u,v) = \{xu + yv / x,y \in \mathbb{Z}\} \, \textit{ and } \, \Lambda' (u',v') = \{xu' + yv' / x,y \in \mathbb{Z}\}$$ 
to be two non-equivalent lattices associated to the same doubly periodic weaving diagram $D_{\infty}$ and let $k_u$ and $k_v$ be two positive integers satisfying, 
$$u' = k_u \times u \, \textit{ and } \, v' = k_v \times v.$$
It follows that $\Lambda$ is contained in $\Lambda'$ ($\Lambda \subseteq \Lambda'$).
Note now that any two weaving motifs $D_W = D_{\infty} \setminus \Lambda$ and $D'_W = D_{\infty} \setminus \Lambda'$ are not necessarily equivalent by Proposition~\ref{R-theorem}, althought they both lift to the same periodic weaving diagram $D_{W, \infty}$ in the plane. This leads to the definition of the notion of \textit{scale-equivalence} between weaving motifs. Note that this inclusion relation guarantees the existence of a \textit{minimal lattice}, denoted by $\Lambda_{min} = \Lambda (u_{min},v_{min})$, such that for any periodic lattice $\Lambda (u,v)$ of $W$, $\Lambda_{min} \subseteq \Lambda (u,v)$.

\begin{definition}\label{def:scale-equivalence}
Let $D_{\infty}$ be a doubly periodic weaving diagram and let $\Lambda_1$ and $\Lambda_2$ be two non-equivalent point lattices such that  $\Lambda_1 \subset \Lambda_2$.
Moreover, let $D_{W_1}= D_{\infty} \setminus \Lambda_1$ and $D_{W_2}= D_{\infty} \setminus \Lambda_2$ be two weaving motifs of $D_{\infty}$. 
Then, $D_{W_1}$ and $D_{W_2}$ are said to be \textit{scale-equivalent} if there exists a weaving motif $D'_{W_1}$ defined as \textit{adjacent} copies of $D_{W_1}$ such that $D'_{W_1} = D_{\infty} \setminus \Lambda_2$ is a weaving motif of $D_{W, \infty}$ for $\Lambda_2$. 
\end{definition}

From now on, we will consider equivalence of doubly periodic weaves by including to scale-equivalence.

\begin{theorem}\label{th:Doubly-R-Th}
Let $W_1$ and $W_2$ be two doubly periodic weaves with weaving diagrams $D_{\infty,1}$ and $D_{\infty,2}$. Let also $\Lambda_1 \subset \Lambda_2$ be two non-equivalent point lattices such that for $i \in \{1,2\}$,  $D_{W_i} = D_{\infty,i} \setminus \Lambda_i$ is a weaving motif of $W_i$. 
Then, $W_1$ and $W_2$ are equivalent if the two following conditions are satisfied, 
\begin{itemize}
    \item $D_{W_1}$ and $D_{W_2}$ are scale-equivalent with $D_{W_1} \subset D'_{W_1}$ and $ D'_{W_1} = D_{\infty,1} \setminus \Lambda_2$, 
    \item $D'_{W_1}$ and $D_{W_2}$ are related by a finite sequence of Reidemeister moves $R_1$, $R_2$, and $R_3$, torus isotopies and Dehn twists. \\
\end{itemize}
\end{theorem}

\subsection{Generalization to hyperbolic weaves and diagrams on higher genus surfaces}\label{sec:2-3}

~

~

\begin{figure}[ht]
\centering
   \includegraphics[width=5.2in]{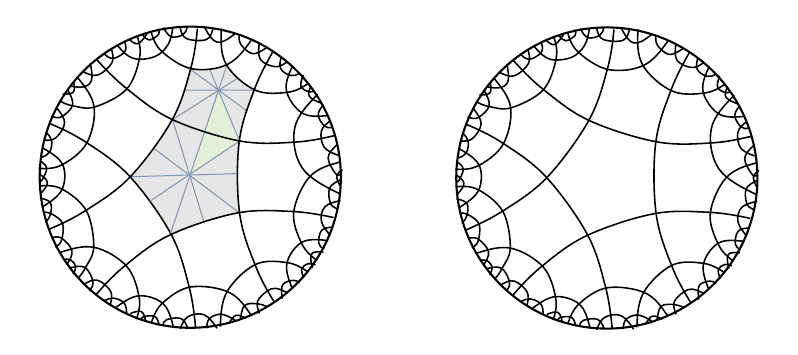}
      \caption{\label{hyperbolic.tiling} Example of kaleidoscopic tiling of the hyperbolic plane.}
\end{figure}

The study of weaving motifs on the torus encourages a generalization to \textit{higher genus surfaces}. This implies a definition of periodic weaving diagrams on the \textit{hyperbolic plane} $\mathbb{H}^2$. However, note that the notion of parallel directions considered to define our sets of threads in $\mathbb{E}^2$ cannot be extended to $\mathbb{H}^2$. We thus restrict our definition of hyperbolic weaves to the cases where their regular projections are isotopic to quadrivalent \textit{kaleidoscopic tilings} of $\mathbb{H}^2$. A kaleidoscopic tiling is constructed by applying repetitively reflection symmetries along the sides of a given hyperbolic convex polygon (Figure~\ref{hyperbolic.tiling}). We will use the Poincare disk representation of $\mathbb{H}^2$ and refer to \cite{Kaleidoscopic} for more details on hyperbolic kaleidoscopic tilings.

\begin{definition}
    Let $P$ be the generating convex polygon of a quadrivalent kaleidoscopic tiling $\mathcal{T}$ of $\mathbb{H}^2$ and let $N \geq 3$ be the number of edges of $P$. 
    If $N$ is odd, then each edge is assigned a different axis of direction. 
    Otherwise, each pair of opposite edges of $P$ are given the same direction.
\end{definition}

\begin{figure}[ht]
\centering
   \includegraphics[width=5in]{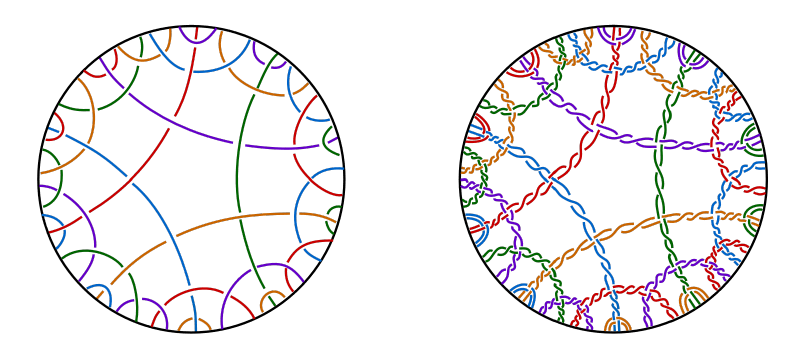}
      \caption{\label{hyperbolic.weave} Example of hyperbolic weaving diagrams.}
\end{figure}

\begin{remark} 
Since reflection symmetries preserves the axis of direction by construction, the notion of sets of threads in $\mathbb{H}^2$ is defined in a coherent way that generalizes the definition in $\mathbb{E}^2$, as shown in Figure~\ref{hyperbolic.weave}.
\end{remark}

Therefore, by specifying each vertex of a hyperbolic kaleidoscopic tilings with an over or under crossing information, we introduce a new class of periodic weaves in $\mathbb{H}^2 \times I$. Note that most of the definitions stated in \cite{Sonia1} may follow naturally here. Two examples of periodic hyperbolic weaving diagrams are illustrated in Figure~\ref{hyperbolic.weave}.

\begin{figure}[ht]
\centering
   \includegraphics[width=5in]{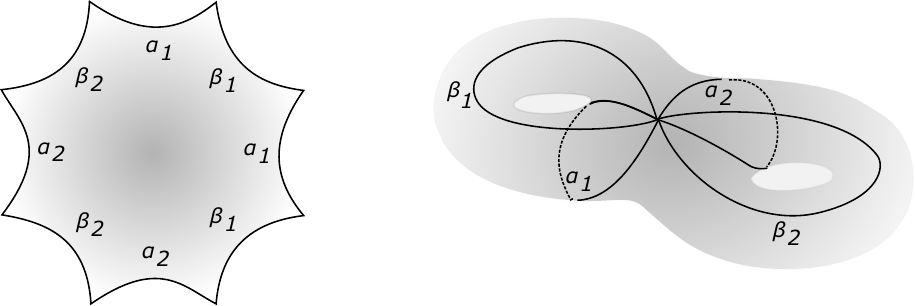}
      \caption{\label{hyperbolic.polygon} A regular $S_2$-tile and its corresponding marked hyperbolic surface.}
\end{figure}

It follows that a hyperbolic weaving motif can be defined as a link diagram on a surface $S_g$ of genus $g \geq 2$. In particular, this results from the pairwise identification of the sides of a generating cell, defined as the quotient of a periodic weaving diagram of $\mathbb{H}^2$ by a discrete lattice. More precisely, a flat weaving motif is a set of simple open arcs and crossings on a hyperbolic polygon. Note that this polygon can be chosen in infinitely many ways, the easiest being a regular $4g$-gons of $\mathbb{H}^2$ (Figure~\ref{hyperbolic.polygon}). This leads to a generalization of the theory of weaving motifs on a torus to motifs on higher genus surfaces. We consider any closed orientable surface $S_g$ of genus $g \geq 2$ with generating loops $\alpha_1, \cdots, \alpha_g$ and $\beta_1, \cdots, \beta_g$. These loops define homotopy classes through a common base point on $S_g$ and map to a basis of $\mathbb{H}^2$.

\begin{definition}\label{def:HyperbolicWeave}
    Let $W$ be a link embedded in $S_g \times I$ of genus $g \geq 2$ with generating loops $\alpha_1, \cdots, \alpha_g$ and $\beta_1, \cdots, \beta_g$, satisfying the two following conditions,
    \begin{enumerate}
        \item each of its components is an essential closed curve that lifts to an unknotted thread,
        \item its components can be partitioned into $N \geq 2$ sets of threads. 
    \end{enumerate}
    Then, the lift of $W$ under the covering map $\rho': \mathbb{H}^2 \times I \rightarrow{} S_g \times I$ is called a \textit{hyperbolic periodic weave} $W_{\infty}$ with $N$ sets of threads. 
    Moreover, the projection of $W$ onto $S_g$ is called a (hyperbolic) \textit{weaving motif}, denoted by $D_W$, and the lift $D_{\infty}$ of $D_W$ under $\rho'$ is called a (hyperbolic) \textit{weaving diagram}.
\end{definition}

Now, to generalize the notion of equivalent weaving motifs on the torus to higher genus surfaces, one needs to consider the \textit{Teichm\"uller space} of a surface $S_g$ of genus $g \geq 2$, and its \textit{mapping class group}. We refer to \cite{Farb} for a complete study of mapping class groups. First, start from any hyperbolic motif, whose flat boundary is a geodesic hyperbolic $4g$-gon on $\mathbb{H}^2$, meaning a polygon such that the sum of its interior angles is equal to $2\pi$. Then, label its edges such that they can be identified pairwise, which results in a closed marked hyperbolic surface $S_g$ of genus $g \geq 2$, as illustrated in Figure~\ref{hyperbolic.polygon}. Such a polygon is called a $S_g$-\textit{tile} and the Teichm\"uller space of the corresponding surface $S_g$ can be seen as the space of marked surfaces homeomorphic to it. Moreover, it is well-known that the Teichm\"uller space of  $S_g$ is in bijection with the set of equivalence classes of hyperbolic $S_g$-tiles. Note that two $S_g$-tiles are said to be \textit{equivalent} if they differ by a marked, orientation-preserving isometry and by \say{pushing the basepoint}, which is the point on the surface where all the vertices of a $S_g$-tile meet after gluing.
The details of this bijection are given in \cite{Farb}. 

Now, to prove the existence and relation between infinitely many $S_g$-tiles, we start from an arbitrary point in $\mathrm{Teich}(S_g)$, the Teichm\"uller space of $S_g$, which represents the equivalence class of a marked surface $S_g$ of genus $g$. From the above bijection, $\mathrm{Teich}(S_g)$ corresponds to the equivalence class of isometric $S_g$-tiles, each corresponding to a unique point in $\mathrm{Teich}(S_g)$. Thus, if there exists a non-isometric $S_g$-tile that can be taken as a unit cell of the same hyperbolic weaving diagram, then it corresponds to a different marked surface and is not isometric to the original $S_g$-tile. However, there exists a relation between these different markings of $S_g$, characterized by its mapping class group $MCG(S_g)$, whose simplest infinite-order elements are known to be \textit{Dehn twists}, as in the torus case. Indeed, given any marked surface $S_g$ of genus $g$ in $\mathrm{Teich}(S_g)$, the marking can be changed by the action of any finite sequence of Dehn twists, which generates non-isometric corresponding $S_g$-tiles. Thus, this proves the existence of infinitely many weaving motifs for any periodic hyperbolic weaving diagram on $\mathbb{H}^2$. In other words, for any given periodic weaving diagram and a corresponding point lattice, every weaving motifs can be obtained from an arbitrarily chosen one by a finite sequence of Dehn twists of $S_g$ along its generating loop, which are extensions of the meridians and longitudes of a torus \cite{Lavasani}. Therefore, to define equivalence class of periodic weaves in terms of \textit{ambient isotopy} in the thickened hyperbolic plane, we can use the above arguments to encode the periodicity and the relation between all possible unit cells for a given lattice, which generalizes Definition~\ref{th:Doubly-R-Th}.

\begin{theorem}\label{th:Periodic-R-Th}
Let $W_{\infty,1}$ and $W_{\infty,2}$ be two periodic hyperbolic weaves with weaving diagrams $D_{\infty,1}$ and $D_{\infty,2}$ in $\mathbb{H}^2$ and let $\Lambda_1 \subset \Lambda_2$ be two non-equivalent point lattices such that for $i \in \{1,2\}$,  $D_{W_i} = D_{\infty,i} \setminus \Lambda_i$ is a weaving motif of $D_{\infty,i}$ on $S_g$. 
Then, $W_{\infty,1}$ and $W_{\infty,2}$ are equivalent if the two following conditions are satisfied, 
\begin{itemize}
    \item $D_{W_1}$ and $D_{W_2}$ are scale-equivalent with $D_{W_1} \subset D'_{W_1}$ and $ D'_{W_1} = W_{01} \setminus \Lambda_2$, 
    \item $D'_{W_1}$ and $D_{W_2}$ are related by a finite sequence of Reidemeister moves $R_1$, $R_2$, and $R_3$, isotopies and Dehn twists of $S_g$. \\
\end{itemize}
\end{theorem}

\subsection{Some particular weaving diagrams}\label{sec:2-4}

~

~

\begin{figure}[ht]
\centering
   \includegraphics[width=5.5in]{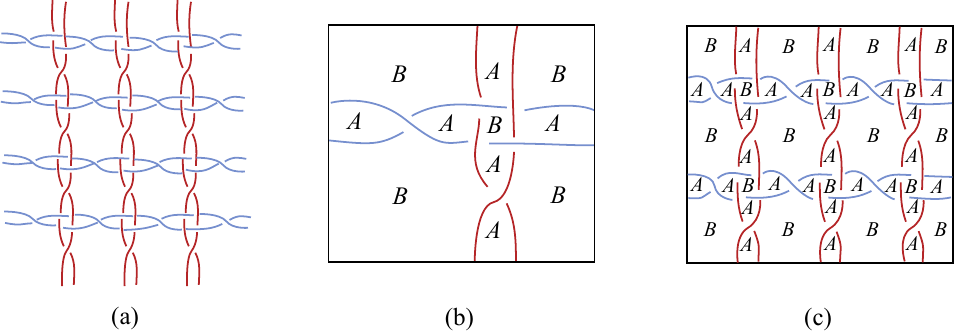}
      \caption{\label{reduced} A periodic reduced (without nugatory crossing) planar diagram (a), a reduced improper weaving motif (b), and a reduced proper weaving motif (c).}
\end{figure}

Many definitions from classical knot theory \cite{Kauffman, MurasugiBook} can be naturally extended for weaving diagrams. In particular, a periodic weave, or any associated weaving diagram or motif, is said to be \textit{alternating} if its crossings alternate cyclically between undercrossings and overcrossings, as one travels along each of its components (see Figure~\ref{reduced}). However, the notion of \textit{reduced} diagrams does not follow directly from links in the thickened surfaces but takes into account the universal cover. More specifically, we have the following. 

\begin{definition}\label{def:2-8}
A weaving motif $D_W$ in $S_g$, with $g \geq 1$, is said to be \textit{reduced} if its lift to $\mathbb{X}^2 = \mathbb{E}^2$ or $\mathbb{H}^2$ does not contain a nugatory crossing. A \textit{nugatory crossing} is a crossing in the diagram so that two of the four local regions at the crossing are part of the same region in the associated infinite diagram. 
\end{definition}

Moreover, at the level of the torus, any crossing $c$ of a weaving motif $D_W$ is called \textit{proper} if the four regions around $c$ delimited by the projection of the threads are all distinct. When every crossing of $D_W$ is proper, $D_W$ is said to be \textit{proper} (see Figure~\ref{reduced}).

Note that an equivalent definition of reduced diagram can be stated using the point lattice inclusion defined above. 

\begin{definition}\label{def:reduced}
    A weaving motif $D_W$ for a fixed point lattice $\Lambda$ is said to be \textit{reduced} if one of the two following conditions is satisfied, 
\begin{itemize}
    \item all its crossings are proper, 
    \item for each improper crossing $c$, there exists a point lattice $\Lambda'$ such that $\Lambda \subset \Lambda'$ and for which $c$ is proper in $D'_W$, where $D'_W$ is the weaving motif constructed by gluing a copy of $D_W$ to each of its boundary sides by translation.
\end{itemize}
\end{definition}

\begin{remark}
    Recently, in \cite{Boden2022}, a different definition of reduced diagram for classic links on a thickened surface has been stated and could also apply here.
\end{remark}


\section{The Bracket Polynomial of Periodic Weaves}\label{sec:3}

\subsection{A Kauffman-type weaving invariant}\label{sec:3-1}

~

~

This section recalls results from \cite{GrishanovP2} and \cite{Grishanov1}, and extends the definition of the bracket polynomial of a weaving motif on a torus to any surface $S_g$ of genus $g \geq 1$.

\begin{definition}\label{def:skein}
Let $D_W$ be a weaving motif on a surface $S_g$ of genus $g \geq 1$, and let $\langle D_W \rangle$  be the element of the ring $\mathbb{Z}[A, B, d]$ defined recursively by the following identities,
  \begin{enumerate} 
  \item $\langle O \rangle = 1$, with $O$ a null-homotopic simple closed curve on $D_W$.
  \item $\langle D_W \cup O \rangle = d \langle D_W \rangle$, when adding an isolated circle $O$ to a diagram $D_W$.
  \item $\langle$ \raisebox{-3pt}{\includegraphics[scale=0.45]{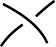}}  $\rangle$ = $A$$\langle$  \raisebox{-0pt}{\includegraphics[scale=0.45]{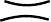}} $\rangle$$+$ $B$$\langle$  \raisebox{-2pt}{\includegraphics[scale=0.55]{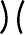}} $\rangle$, for diagrams that differ locally around a single crossing.  
  \end{enumerate}
This last relation is called the \textit{skein relation} and $\langle D_W \rangle$ denotes the \textit{bracket polynomial}.
\end{definition}

\begin{figure}[ht]
\centering
   \includegraphics[width=5.5in]{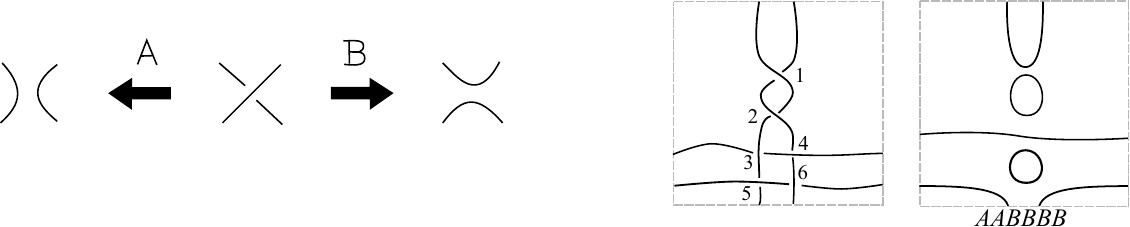}
      \caption{\label{state} On the left, the two types of splitting.  On the right, an example of a state $S=AABBBB$ of a weaving motif.}
\end{figure}

This polynomial is known to be well defined on classic link diagrams and can be generalized to periodic weaves. Let $D_W$ be any weaving motif of a periodic weave of $\mathbb{X}^2 = \mathbb{E}^2$ or $\mathbb{H}^2$. Then, every crossing of $D_W$ can be smoothed via an operation of type $A$ or $B$, as illustrated in Figure~\ref{state}. The overall operations can be expressed as a state $S$ of $D_W$, defined as a sequence of symbols $A$ and $B$ of length $C$, where $C$ is the number of crossings of $D_W$, 
$$S = A B A A B B … A B B A.$$

It is well-known for links in $\mathbb{S}^3$ that a diagram $D_S$ in a state $S$ is a disjoint union of $c_S$ null-homotopic simple closed curves. It follows from Definition~\ref{def:skein}(2) that,  
$$\langle D_S \rangle = d\,^{c_{S-1}}.$$
Moreover, if $i$ is the number of splits of type $A$ and $j$ the number of splits of type $B$, then the total contribution of the state $S$ to the bracket polynomial is given by applying the skein relation recursively,
$$P_S = \langle D_S / S \rangle  \langle D_S \rangle = A^i B\,^j d\,^{c_{S-1}}.$$ 

However, recall that a state of a weaving motif $D_W$ may contain essential simple closed curves, as in the case of any link diagram on a surface $S_g$ of genus $g \geq 1$. Such a set is called a \textit{winding} in \cite{GrishanovP2} and is denoted by $(m,n)_S$ for $g=1$, where $m$ and $n$ are the number of intersections of the winding with a torus meridian and longitude, respectively. For example, in the case of the state diagram of Figure~\ref{state} (right), $(AABBBB)$, we have $(m, n)_S = (0, 2)_S$. 

For the general case of $g \geq 1$, a winding is denoted by, 
$$(m_1, \cdots, m_g, n_1, \cdots, n_g)_S ,$$ 
where $m_1, \cdots, m_g, n_1, \cdots, n_g$ are the number of intersections of the winding with the generating loops $\alpha_1, \cdots, \alpha_g, \beta_1, \cdots, \beta_g$ of the surface $S_g$ respectively, see Figure~\ref{hyperbolic.polygon} and \cite{Lavasani} for more details.

Recall that by definition, windings also encode the periodicity of weaves in the universal cover, which is not considered for classic links in a thickened surface, as in \cite{AdamsFleming, Boden2019, Boden2022, Fleming}. We can thus generalize the value of the bracket polynomial of a winding defined for $g=1$ in \cite{GrishanovP2}, with respect to Definition~\ref{def:skein},  
$$\langle (m_1, \cdots, m_g, n_1, \cdots, n_g)_S \rangle \ = \  (m_1, \cdots, m_g, n_1, \cdots, n_g)_S, \mbox{ for every } g \geq 1.$$

Therefore, following the above reasoning, the bracket polynomial of a weaving motif is well-defined.

\begin{proposition}\label{prop:3-2}
The bracket polynomial $\langle D_W \rangle$ of a weaving motif $D_W$ on a surface $S_g$ of genus $g\geq 1$ is uniquely determined by the identities (1), (2), (3) of Definition~\ref{def:skein}, and is expressed by summation over all states of the diagram,  
\begin{equation}                        
\langle D_W \rangle = \sum_{S} A^i B\,^j d\,^{c_S} (m_1, \cdots, m_g, n_1, \cdots, n_g)_S. 
\end{equation}
\end{proposition}

As for links in $\mathbb{S}^3$, the bracket polynomial of a weaving motif is proven to be invariant under the Reidemeister moves following the same approach of \textit{Lemma 2.3} in \cite{Kauffman}.

\begin{lemma}\label{lem:3-3}
If the three diagrams represent the same weaving motif except in the area indicated, we have 
$\langle$ \raisebox{-2pt}{\includegraphics[scale=0.45]{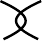}}  $\rangle$ = $AB$$\langle$  \raisebox{-2pt}{\includegraphics[scale=0.55]{skein_3.pdf}} $\rangle$ $+$ $(ABd + A^2 + B^2)$$\langle$  \raisebox{-0pt}{\includegraphics[scale=0.5]{skein_2.pdf}} $\rangle$.  
\end{lemma}

Thus, the bracket is invariant for the Reidemeister move $R_2$ for all diagrams if
$$AB = 1 \mbox{ and } d = -A^2 - A^{-2}.$$

Moreover, this implies also the invariance of the bracket for the Reidemeister move $R_3$, which concludes the invariance by \textit{regular isotopy}.

\begin{lemma}\label{lem:3-4}
The bracket invariance for the Reidemeister move $R_2$ implies the bracket invariance for the Reidemeister move $R_3$. 
Thus, the bracket polynomial is an invariant of regular isotopy for periodic weaves for a fixed point lattice.
\end{lemma}

Finally, to prove the invariance uner the Reidemeister move $R_1$, we use the following proposition that provides an identity for $R_1$, as in \cite{Kauffman}.

\begin{proposition}\label{prop:3-5}
If $AB = 1$ and $d = -A^2-A^{-2}$, then, for the Reidemeister move $R_1$, we have
$$
 \left\{
    \begin{array}{ll}
        \langle \ \raisebox{-3pt}{\includegraphics[scale=0.4]{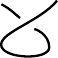}} \ \rangle = (- A^3) \langle \ \raisebox{1pt}{\includegraphics[scale=0.4]{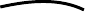}} \ \rangle , \\
        ~ \\
        \langle \ \raisebox{-3pt}{\includegraphics[scale=0.38]{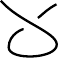}} \ \rangle = (- A^{-3}) \langle \ \raisebox{1pt}{\includegraphics[scale=0.4]{353.pdf}} \ \rangle . \\
    \end{array}
\right.
$$
~
\end{proposition}
 
\begin{figure}[ht]
\centering
   \includegraphics[scale=0.50]{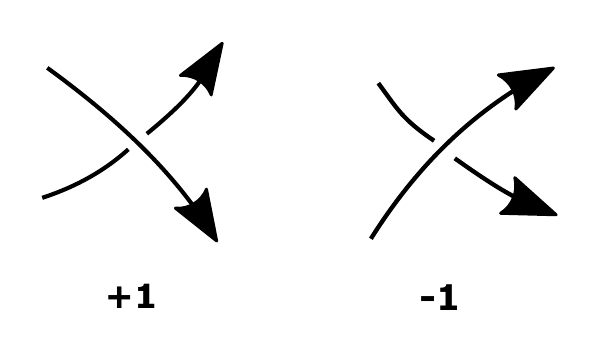}
      \caption{\label{sign} Sign convention.}
\end{figure}

The strategy is to use the \textit{writhe} $w_r(D_W)$ of a weaving motif $D_W$, which is the sum of the signs 
of all the crossings, where each crossing is given a sign $\pm 1$, as in Figure~\ref{sign}. 

Any weaving motif $D_W$ consists of $T$ essential closed curve components, each denoted by $t_i$, that can be oriented in an arbitrary way. 
We call $D^i_W$ the part of the diagram $D_W$ that corresponds to the component $t_i$. 
Then we have in $D_W$, 
\begin{equation}
w_r(D_W) = \sum_{i=1}^{T} w_r(D^i_W).
\end{equation}
~ \\

We can now define a polynomial constructed from the bracket. 
For every $g \geq 1$, we set
\begin{equation}
\begin{aligned}
f(D_W) &= (-A)^{-3w_r(D_W)} \langle D_W \rangle , \\
&= (-A)^{-3w_r(D_W)} \Bigl(\sum_{S} A^{i-j} (-A^2 – A^{-2})^{c_S} (m_1, \cdots, m_g, n_1, \cdots, n_g)_S \Bigr).\\
\end{aligned}
\end{equation}
~ \\

\begin{theorem}\label{th:3-6}
The polynomial $f(D_W) \in \mathbb{Z}[A]$ defined above is an ambient isotopic invariant for oriented perioding weaves under a fixed point lattice.
\end{theorem}

\begin{proof}
From Lemma~\ref{lem:3-4}, we already have the invariance of $f(D_W)$ for the Reidemeister moves $R_2$ and $R_3$. Then, by combining the behavior of the writhe defined above under the Reidemeister move $R_1$ with the previous relation of the bracket for $R_1$ in Proposition~\ref{prop:3-5}, it follows that $f(D_W)$ is invariant under $R_1$ type moves. Thus, $f(D_W)$ is invariant under all three moves, and is therefore an invariant of ambient isotopy for the chosen lattice.
\end{proof}

Nevertheless, this polynomial still depends on the choice of the unit cell for the chosen point lattice, since the multipliers $(m_1, \cdots, m_g, n_1, \cdots, n_g)_S$ describing the windings depend on the Dehn twists of the surface $S_g$. As seen earlier, to have a weaving invariant for a fixed lattice, we also need the invariance of the polynomial under Dehn twists of $S_g$. Once again, the particular case of the torus is described in \cite{GrishanovP2} and is generalized below to $g \geq 1$.

\begin{theorem}\label{th:3-7}
The polynomial $f(D_W)$, when $(m_1, \cdots, m_g, n_ 1, \cdots, n_g)_S$ is in a \textit{canonical form} for each state $S$, defines a Kauffman-type weaving invariant for a fixed point lattice.
\end{theorem}

\begin{proof}
To construct an invariant independent of the Dehn twists of $S_g$, a possibility is to define a canonical form for the set $\{v_S\}$ $=$ $\bigl\{ (m_1, \cdots, m_g, n_1, \cdots, n_g)_S \bigr\} $ of windings for every state $S$.  
Indeed, since this set $\{v_S\}$ depends on the twists of $S_g$, one must transform it into the canonical form to make it invariant. The \textit{Dehn-Lickorish Theorem} states that it is sufficient to select a finite number of Dehn twists to generate the mapping class group $MCG(S_g)$ of a surface $S_g$ of genus $g$. Moreover, since the map $\psi: MCG(S_g) \to Sp(2g, \mathbb{Z})$ is surjective for $g \geq 1$, it follows that the images of the Dehn twists generate $Sp(2g, \mathbb{Z})$ (\cite{Farb}). Besides, recall that the determinant of every matrix $U \in Sp(2g, \mathbb{Z})$ is equal to $1$ and that for $g = 1$, $Sp (2g, \mathbb{Z}) = SL_2(\mathbb{Z})$. Thus, following the strategy used in \cite{GrishanovP2}, one can represent the transformation of a winding $v_S = ((m_1, \cdots, m_g, n_1, \cdots, n_g)_S$ by a sequence of Dehn twists of $S_g$ as a product of $v_S$ by a matrix $U \in Sp (2g, \mathbb{Z})$, $$v'_S = v_S . U ,$$ considering the canonical matrix multiplication on $Sp(2g, \mathbb{Z})$.

To define the canonical form of a set $V = \{v_S\}$, we associate a quadratic functional $Q$, 
\begin{equation}
Q(V) := \sum_{S}^{N} |v_S|^2.
\end{equation}
Thus, a finite sequence of Dehn twists, given by a matrix $U$, converts the set of windings $V = \{v_S\}$ to a set $V' = \{v'_S\}$, with $v'_S = v_S . U$ and the value of $Q$ becomes,         
\begin{equation}                
Q(V') = \sum_{s}^{N} |v'_S|^2 = \sum_{S}^{N} v_S . U . U^T . v_S^T.
\end{equation}
Then, to construct a canonical form of the windings $V = \{v_S\}$, the idea is to find a finite sequence of Dehn twists encoded in $U$ that minimizes the value of $Q$ (\cite{GrishanovP2}),
\begin{equation}
Q(V') = \sum_{S}^{N} |v'_S|^2 = \sum_{S}^{N} v_S . U . U^T . v_S^T \longrightarrow \min, U \in Sp(2g, \mathbb{Z}).
\end{equation}
This equation always has a unique non-trivial solution $U_0$. Indeed, let $M = U . U^T$ be a symmetric definite positive matrix. Then, if $x = v_S$ and $\phi(x) = x . M . x^T$, then $\phi(0)=0$ and for all $x \neq 0$, $\phi(x) > 0$.
Thus, there exists an orthonormal basis $\{e_1, \cdots, e_d\}$ such that for all $i$ in  $\{1,\cdots ,d \}$, $e_i$ is an eigenvector of $M$. We denote the corresponding eigenvalue $\lambda_i$ and we show that $\phi$ is strictly convex.

Let $0 < \mu < 1$ and consider $\phi \bigl(\mu x + (1-\mu)y\bigr)$, with $x \neq y$. Then,

\begin{equation}
\begin{aligned}
\phi\bigl(\mu x + (1-\mu)y\bigr) &= \langle \mu x + (1-\mu)y, M (\mu x + (1-\mu)y) \rangle , \\
&= \langle \sum_{i=1}^{d} \mu x_i e_i + (1-\mu)y_i e_i, M \sum_{i=1}^{d} \mu x_i e_i + (1-\mu)y_i e_i \rangle ,\\
&= \sum_{i=1}^{d} \lambda_i \bigl(\mu x_i + (1-\mu)y_i\bigr)^2.
\end{aligned}
\end{equation}
Moreover, $x^2$ is strictly convex and for some $i$, $x_i \neq y_i$, thus,

\begin{equation}
\sum_{i=1}^{d} \lambda_i \bigl(\mu x_i + (1-\mu)y_i\bigr)^2 < \sum_{i=1}^{d} \lambda_i \bigl(\mu x_i^2 + (1-\mu)y_i^2\bigr).
\end{equation}
Therefore, since $\phi$ is strictly convex and has a limit at infinity, it has a unique minimum, which concludes our proof. So, for every state $S$, the canonical form of a set of windings $V = \{v_S\}$, with the winding as coordinates 
$v_S =(m_1, \cdots, m_g, n_ 1, \cdots, n_g)_S$ is an invariant and thus, $f(D_W)$ too.
\end{proof}
~ \\

\subsection{The case of alternating weaving diagrams}\label{sec:3-2}

~

~

Now, we study the bracket polynomial for the case of alternating weaves. It is well-known that the \textit{degree} of a polynomial is the most important aspect of the polynomial as an invariant \cite{Murasugi2}. The following proposition and its proof follow the strategy of the similar result for classic links in $\mathbb{S}^3$, but depend strongly on the definition of reduced and proper diagrams stated in Section~\ref{sec:2-4}.

\begin{proposition}\label{prop:3-8}
Let $D_W$ be an alternating reduced weaving motif colored so that all the regions labeled $A$ are white and all the regions labeled $B$ are black. Let $\mathit{C}$ be the number of crossings, $W$ be the number of white regions and $B$ be the number of black ones. Then,              
\begin{equation}
\begin{aligned}
\max \mathrm{deg} (\langle D_W \rangle) &= \mathit{C} + 2W,\\
\min \mathrm{deg}  (\langle D_W \rangle) \, &= - \mathit{C} - 2B,
\end{aligned}
\end{equation}
with $\max \mathrm{deg} (P)$ and $\min \mathrm{deg} (P)$ are respectively the maximal and the minimal degree of any polynomial $P$ in $\mathbb{Z}[A, B, d]$.  
\end{proposition}

\begin{proof}
Since $D_W$ is alternating, it admits a canonical checkerboard coloring by definition, which means that two edge-adjacent regions always have different colors. Let $S = S_A$ be the state obtained by splitting every crossing in the diagram $D_W$ in the $A$-direction. Then we have $\langle D_W / S \rangle = A^{\mathit{C}}$, and since the number of components $c_S$ is equal to $W$, thus as seen earlier, the total contribution of the state $S$ to the bracket polynomial is given by, 
\begin{equation}
\begin{aligned}
P_S &= \langle D_W / S \rangle  d\,^{c_S} (m_1, \cdots, m_g, n_1, \cdots, n_g)_S \\
&= A^{\mathit{C}}  d\,^W (m_1, \cdots, m_g, n_1, \cdots, n_g)_S,\ g \geq 1.
\end{aligned}
\end{equation}
And since $d = -A^2-A^{-2}$, then $\max \mathrm{deg} (P_S) = \mathit{C}+2W$, which is the desired relation.  

Now let $S’ \ne S_A$ be any other state and verify that $\mathrm{deg} (P_{S’}) \ngeq \mathrm{deg} (P_S)$. Then $S’$ can be obtained from $S = S_A$ by switching some splittings of $S$. A sequence of states can be defined by $S(0), S(1), \cdots, S(n)$ such that $S = S(0)$, $S’ = S(1)$. Thus, for every positive integer $i$, $S(i+1)$ is obtained from $S(i)$ by switching one splitting from type $A$ to type $B$, which then contributes a factor of $A^{-1}$ in the polynomial
\begin{equation}
\langle D_W /S(i+1) \rangle = A^{-2} \langle D_W /S(i) \rangle.
\end{equation}

We now need to distinguish two cases.

Case 1. The weaving motif $D_W$ is reduced and proper for the fixed point lattice. Then, $c_{S(i+1)} \leq c_{S(i)}-1$, since switching one splitting can change the component number by at most one. Thus, $\max \mathrm{deg} (P_{S(i+1)}) \leq \max \mathrm{deg} (P_{S(i)})$. Moreover, let $c$ be the crossing point for which we change the $A$-splice into the $B$-splice from $S(0)$ to $S(1)$. Since $D_W$ is proper, the crossing $c$ is proper. Thus, we can use the following lemma, (\textit{Lemma 3.2} in \cite{Kamada}).

\begin{lemma}\label{lem:3-9}
Let $D_W$ be an alternating weaving motif and let $S_A$ (resp. $S_B$) be the state of $D_W$ obtained from $D_W$ by doing an $A$-splice (resp. $B$-splice) for every crossing. For a crossing $c$ of $D_W$, let $R_1(c)$ and $R_2(c)$ be the closed regions of $S_A$ (or $R’_1(c)$ and $R’_2(c)$ be the closed regions of $S_B$) around $c$. If $c$ is a proper crossing, then $$R_1(c) \ne R_2(c) \textit{ and } R’_1(c) \ne R’_2(c).$$
 \end{lemma}

\begin{proof}
Since $c$ is a proper crossing, the four closed regions of $D_W$ appearing around $c$ are all distinct. Moreover, since $D_W$ is alternating, it has a canonical checkerboard coloring and there is a one-to-one correspondence,
$$\bigl\{ \mbox{the closed regions of $S_A$} \bigr\} \cup \bigl\{ \mbox{the closed regions of $S_B$} \bigr\} \to \bigl\{ \mbox{the closed regions of $D_W$} \bigr\} $$
Then $R_1(c), R_2(c), R’_1(c) \textit{ and } R’_2(c)$ correspond to the four distinct closed regions of $D_W$ around $c$. This concludes the proof.
\end{proof}

Thus, from this lemma, since $S(1)$ is obtained from $S(0)$ by changing an $A$-splice to a $B$-splice at $c$, two distinct regions $R_1(c)$ and $R_2(c)$ become a single region. 
Hence $c_{S(1)} = c_S - 1$. To conclude, the term of maximal degree in the entire bracket polynomial is contributed by the state $S = S_A$, and is not canceled by terms from any other state, so we arrive at
$$\max \mathrm{deg} (\langle D_W \rangle) = \mathit{C} + 2W. $$ 
The proof is similar for,
$$\min \mathrm{deg} (\langle D_W \rangle) \, = - \mathit{C} - 2B.$$

Case 2: The weaving motif is reduced but not proper for the fixed lattice. Then, there exists at least one crossing which is not proper. If we change an $A$-splice to a $B$-splice at a crossing $c$ that is proper, then the conclusion is the same as before. Now, if we change an $A$-splice to a $B$-splice at a crossing $c’$ that is not proper, then some white regions would touch both sides of a crossing. In this case, the number of split components does not decrease from $S(0)$ to $S(1)$, $c_{S(1)} = c_{S(0)}$. But, as seen before, $\langle D_W /S(1) \rangle = A^{-2} \langle D_W /S(0) \rangle$ and there is no isthmus in the diagram, so the number of components either decreases or is constant. Therefore, 
\begin{equation}
\max \mathrm{deg} (P_{S(1)}) \leq \max \mathrm{deg} (P_{S(0)}).
\end{equation}
Thus, once again, the term of maximal degree in the entire bracket polynomial is contributed by the state $S = S_A$, and is not canceled by terms from any other state and therefore, 
$$\max \mathrm{deg} (\langle D_W \rangle) = \mathit{C} + 2W. $$ 
The proof is similar for,
$$\min \mathrm{deg} (\langle D_W \rangle) \, = - \mathit{C} - 2B.$$ 
\end{proof}

\begin{figure}[ht]
\centering
   \includegraphics[width=5.5in]{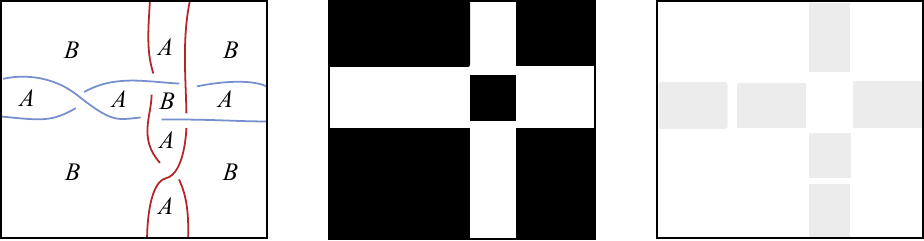}
      \caption{\label{fig10} Example of $D_W$ (left), $S_A = A \cdots A$ (middle) and $S_B = B \cdots B$ (right).}
\end{figure}

Now, it is possible to define a relation between the closed regions of $D_W$ and the regions of the diagram after splitting as done in Section 3 of \cite{Kamada}. Let $S_A$ (resp. $S_B$) be again the state obtained by splitting every crossing in the diagram in the $A$ (resp. $B$)-direction, and $D_W$ be colored so that all the regions labeled $A$ are white (gray on Figure~\ref{fig10}) and all the regions labeled $B$ are black. 

Therefore, we have the following correspondences, \\

$\bigl\{ \mbox{the closed regions of $S_A$} \bigr\}  \longrightarrow \bigl\{ \mbox{the closed regions of $D_W$ in black regions $B$} \bigr\}$,

$\bigl\{ \mbox{the closed regions of $S_B$} \bigr\}  \longrightarrow \bigl\{ \mbox{the closed regions of $D_W$ in white regions $W$} \bigr\}$, \\

which leads to the following bijection, 

$\bigl\{ \mbox{the closed regions of $S_A$} \bigr\} \ \cup \ \bigl\{ \mbox{the closed regions of $S_B$} \bigr\} \longrightarrow \bigl\{ \mbox{the closed regions of $D_W$} \bigr\}.$\\

And when considering a diagram on a surface of genus $g \geq 1$, using the Euler characteristic and the fact that such a diagram is quadrivalent, we conclude that, 

$$\textit{the number of closed regions of } \, D_W \, \textit{ is equal to } \, \mathit{C} +2 - 2g, \, \textit{ for every } \, g \geq 1.$$

It is now possible to extend the proof of \textit{Theorem~2.10} in \cite{Kauffman} to reduced alternating weaving motifs of $\mathbb{X}^2 = \mathbb{E}^2$ or $\mathbb{H}^2$ defined for a fixed integer lattice.

\begin{theorem}\label{th:3-10}
Let $W_{\infty}$ be an alternating periodic weave in the thickened plane $\mathbb{X}^2 \times I$, where $\mathbb{X}^2= \mathbb{E}^2$ or $\mathbb{H}^2$. 
Then, the number of crossings $\mathit{C}$ in an alternating weaving motif $D_W$ is a topological invariant of its corresponding periodic weave $W_{\infty}$ for a fixed point lattice. Therefore any two reduced alternating weaving motifs of a given periodic weave quotient by the same point lattice have the same number of crossings.
\end{theorem}

\begin{proof}
Let $\mathrm{span}(D_W)$ defined by 
$$\mathrm{span}(D_W) = \max \mathrm{deg} (\langle D_W \rangle) - \min \mathrm{deg} (\langle D_W \rangle).$$

Then we have, 
$\mathrm{span}(D_W)=2\mathit{C}+2(W+B)=2\mathit{C}+2(\mathit{C}+2-2g)$.

So finally, $\mathrm{span}(D_W) = 4 \mathit{C} - 4g + 4$, for every $g \geq 1$.
\end{proof}

\begin{remark}
  Note that Theorem~\ref{th:3-10} implies that the number of crossings on a generating cell of a periodic weave is a \textit{topological invariant} for a given point lattice. However, it does not allow a comparison of scale-equivalent weaving motifs, which is an important point to remember for the classification of periodic weaves in a thickened plane since the number of crossing in a unit cell describes the complexity of the structure. \\  
\end{remark}


\section{The Jones Polynomial and Tait’s First Conjecture for weaves}\label{sec:4}

\subsection{The Jones polynomial of a periodic weave}\label{sec:4-1}

~

~

The Jones polynomial is defined by the following identities in \textit{Section 2} of \cite{Kauffman},
 \begin{itemize}
    \item  $J_O = 1$,
    \item  $t\,^{-1} \, J \raisebox{-6pt}{\includegraphics[scale=0.10]{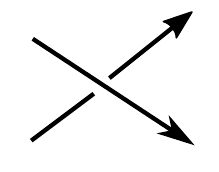}} - t \, J \raisebox{-6pt}{\includegraphics[scale=0.10]{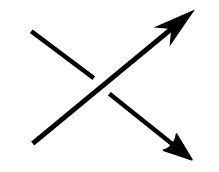}} =  (\sqrt{t} - \frac{1}{\sqrt{t}}) \, J \raisebox{-6pt}{\includegraphics[scale=0.10]{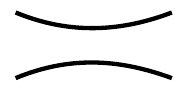}}$.\\
 \end{itemize}

And it is related to the weaving invariant defined above by the following relation, 

\begin{theorem}\label{th:4-1}
The Jones polynomial $J_W$ of a weaving motif $\mathrm{D_W}$ is related to its bracket-type polynomial, for every $g \geq 1$, by the expression, $J_W(t) = f (D_W)(t^{-1/4})$
$$J_{W}(t) = (-1)^{-3w_r(D_W)}  t^{-(1/4)-3w_r(D_W)} \Bigl(\sum_{S} t\,^{i-j} (-t\,^2 - t\,^{-2})^{c_S} (m_1, \cdots, m_g, n_1, \cdots, n_g)_S  \Bigr).$$
\end{theorem}

\begin{proof}
By the skein relation:  
$$
 \left\{
    \begin{array}{ll}
        \langle \ \raisebox{-2pt}{\includegraphics[scale=0.8]{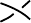}} \ \rangle = A \langle \ \raisebox{0pt}{\includegraphics[scale=0.5]{skein_2.pdf}} \ \rangle + A^{-1} \langle \ \raisebox{-2.5pt}{\includegraphics[scale=0.6]{skein_3.pdf}} \ \rangle, \\
        \langle \ \raisebox{-2pt}{\includegraphics[scale=0.85]{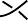}} \ \rangle = A^{-1} \langle \ \raisebox{0pt}{\includegraphics[scale=0.5]{skein_2.pdf}} \ \rangle + A \langle \ \raisebox{-2.5pt}{\includegraphics[scale=0.6]{skein_3.pdf}} \ \rangle . 
    \end{array}
\right.
$$

Thus, we have  $A \langle \ \raisebox{-2pt}{\includegraphics[scale=0.8]{4_1_1.pdf}} \ \rangle - A^{-1} \langle \ \raisebox{-2pt}{\includegraphics[scale=0.8]{4_1_2.pdf}} \ \rangle = (A^2 - A^{-2}) \langle \ \raisebox{0pt}{\includegraphics[scale=0.5]{skein_2.pdf}} \ \rangle $.

If we consider the writhe $w_r(D_W)$ of the weaving motif in the bracket on the right side of the equation, then the other two diagrams on the left have writhes 
$(w_r(D_W) + 1)$ and $(w_r(D_W) - 1)$ respectively. Thus, by multiplying the previous equation by the appropriate writhe, we obtain

$$A^4 f(\raisebox{-6pt}{\includegraphics[scale=0.20]{JP1.pdf}}) - A^{-4} f(\raisebox{-6pt}{\includegraphics[scale=0.20]{JP2.pdf}}) 
= (A^{-2} - A^2) f(\raisebox{-2pt}{\includegraphics[scale=0.20]{JP3.pdf}}). $$
\end{proof}
~ \\

\subsection{Tait’s First Conjecture for periodic weaves}\label{sec:4-2}

~

~

Before stating the main result of this paper, it is necessary to give a last essential definition, that is particular to the case of 
weaving motifs of infinite periodic weaving structures. 

\begin{definition}\label{def:4-2}
The \textit{crossing number} $C(W_{\infty})$ of a periodic weave $W_{\infty}$ with weaving diagram $D_{\infty}$ is defined as the minimum number of crossings that can be found among all possible weaving motifs that lift to $W_{\infty}$. In other words, for a minimal point lattice $\Lambda_{min}$,
$$C(W_{\infty}) = \min \bigl\{ C (D_W), \, D_W= W_0 \setminus \Lambda_{min} \bigr\}.$$
Any weaving motif of $W_{\infty}$ which has exactly $\mathit{C(W_{\infty})}$ crossings is said to be a \textit{minimum} motif. 
\end{definition}

It is important to recall at this point that any weaving motif $D_W$ must encode the alternating property and the periodicity of its corresponding weave $W_{\infty}$. 
Moreover, as seen earlier, a minimal diagram of a weave is not unique by construction and it is necessary to identify the minimal point lattice to apply Tait's first and second conjectures to weaving structures.

\begin{theorem}\label{th:4-3} \textbf{(Tait’s First Conjecture for Periodic Alternating Weaves)} 
A minimal reduced alternating weaving motif on a surface $S_g$ of genus $g \geq 1$ is a minimum diagram of its alternating periodic weave in $\mathbb{X}^2 \times I$, where $\mathbb{X}^2= \mathbb{E}^2$ or $\mathbb{H}^2$. 
\end{theorem}

\begin{proof}
Since $J_W(t) = f (D_W)(t^{-1/4})$ and $\mathrm{span}(D_W) = 4\mathit{C} - 4g +4$, for every $g \geq 1$, thus,
\begin{equation}
\begin{aligned}
\mathrm{span}\bigl(J_W(t)\bigr) &= \max \mathrm{deg} \bigl(  J_W(t) \bigr) - \min \mathrm{deg} \bigl( J_W(t) \bigr) .\\
&= \mathit{C} - g - 1.
\end{aligned}
\end{equation}

And the number of crossings is an invariant thus, it is fixed here for a minimal reduced alternating weaving motif. Moreover, we have a generalization of the previous result for the general case, not necessary alternating, that can be proven in a similar way than in the proof of the generalization of \textit{Proposition 2.9} in \cite{Kauffman},
\begin{equation}
\mathrm{span}(D_W) \leq 4 \mathit{C} - 4g + 4, \mbox{ for every } g \geq 1.
\end{equation}

Thus, the number of crossing points cannot decrease below $\mathrm{span}\bigl(J_W(t)\bigr)$.
We conclude that $D_W$ must be a minimum diagram. 
\end{proof}
~ \\


\section{Tait’s Second Conjecture for periodic weaves}\label{sec:5}

~

In this last section, we generalize Tait's second conjecture to periodic alternating weaves in the thickened Euclidean or hyperbolic plane. This result concerns the invariance of the writhe for minimal reduced alternating weaving motifs on a surface $S_g$ of genus $g \geq 1$. The strategy follows the proof of this same conjecture for classic links in $\mathbb{S}^3$ from \cite{Stong}, with a special attention to the fact that the components of a weave are unknotted simple open curves. We will thus only consider the particular case of \textit{self-crossing} of a thread, defined as a finite sequence of Reidemeister moves $R_1$.\\

\subsection{Writhe, linking number and adequacy of weaving diagrams}\label{sec:5-1}

~

~

In the proof of  Theorem.~\ref{th:3-6}, we have seen that the writhe is not invariant under Reidemeister moves of type $R_1$, which only concerns cases of self-crossings in a same thread, as defined above. Therefore, we can start to study crossings between two distinct components of a weaving motif.

\begin{definition}\label{def:5-1}
Let $D_W$ be a weaving motif of an oriented periodic weave $W_{\infty}$ for a fixed point lattice $\Lambda$. Let $t_i$ and $t_j$ be two arbitrary threads of $W_{\infty}$ and denote by $D^i_W$ and $D^j_W$ their image on $D_W$ under the covering map for $\Lambda$. The linking number of $D^i_W$ and $D^j_W$, denoted $lk(D^i_W,D^j_W)$, is the sum taken over crossings between $D^i_W$ and $D^j_W$, where each crossing is assigned a symbol $\pm 1$ according to the convention of Figure~\ref{sign}.
\end{definition}

It is important to specify that the linking number is defined for pairs of threads in a weave. We will now prove that the linking number is a weaving invariant for a fixed point lattice. 

\begin{proposition}\label{prop:5-2}
Let $D_{W_1}$ and $D_{W_2}$ be two weaving motifs defined as the quotient of an oriented periodic weave $W_{\infty}$ by a fixed point lattice $\Lambda$, which differ by a Reidemeister move $R_2$, or $R_3$. Then, with the same notation as before, we have, 
$$lk(D^i_{W_1},D^j_{W_1}) = lk(D^i_{W_2},D^j_{W_2}).$$ 
\end{proposition}

As for classic links, this proof is immediate from the invariance of the writhe under Reidemeister moves of type $R_2$ and $R_3$. Besides, isotopies and Dehn twists on the surface $S_g$ will clearly not affect the linking number by definition. Therefore, due to the above considerations, the linking number is an invariant of regular isotopy and it is thus possible to extend its definition to a thickened surface and to its universal cover for a fixed lattice. 

\begin{definition}\label{def:5-3}
Let $D_W$ be a weaving motif of an oriented periodic weave $W_{\infty}$ for a fixed point lattice $\Lambda$. With the same notation as above, for any two threads $t_i$ and $t_j$ of $W_{\infty}$, the linking number is defined by, 
$$lk(t_i, t_j) := lk(D^i_W,D^j_W)_{\Lambda}.$$
\end{definition}

Then, we can also apply the notion of adequacy of a weaving motif using once again the notion of states described in Section~\ref{sec:3-1} and by following the definition of an adequate link diagram on a surface from \cite{Boden2022}.  

\begin{definition}\label{def:5-4}
Let $D_W$ be a weaving motif of an oriented periodic weave $W_{\infty}$ for a fixed point lattice. Let $S_A$ denote the state of $D_W$ in which all crossings are $A$-smoothed, and $S_B$ the state in which all crossings are $B$-smoothed. For any state $S$ of $D_W$, $c_S$ denote the number of null-homotopic components and $t_S$ the number of components of the winding, if any. If for all states $S$ adjacent to $S_A$, we have $c_{S_A} \geq c_S$ or $t_{S_A} \neq t_S$, then $D_W$ is said to be \textit{A-adequate}. If, for all states $S$ adjacent to $S_B$, we have $c_{S_B} \geq c_S$ or $t_{S_B} \neq t_S$, then $D_W$ is said to be \textit{B-adequate}. If $D_W$ is both A-adequate and B-adequate, then $D_W$ is said to be \textit{adequate}. 
\end{definition}

There exists a method to verify if a weaving motif $D_W$ is A-adequate, that refers to the proof of Proposition.~\ref{prop:3-8} and to \cite{Stong}. Considering the state $S_A$ of $D_W$, we arbitrarily choose a crossing that we switch from an $A$-smoothing to a $B$-smoothing. Such an operation will either decrease or preserve the number of components as detailed in \cite{Boden2022}, unless the chosen crossing was a (positive) self-crossing since it would split a curve into two. Thus, if each component of $S_A$ never forms a (positive) self-crossing at a former crossing of $D_W$, then it is A-adequate. In a similar way, if each component of $S_B$ never forms a (negative) self-crossing at a former crossing of $D_W$, then it is B-adequate. We can deduct immediately from this observation that a reduced weaving motif is always adequate by definition, since it never contains any self-crossing by definition.

One of the key points in the proof of Tait's second conjecture for links in $\mathbb{S}^3$ by R. Stong (\cite{Stong}) was to use the notion of parallels of diagrams and study their adequacy. This strategy has also been used to prove Tait's conjectures for links in thickened surfaces in \cite{Boden2022}. 

\begin{definition}\label{def:5-5}
Let $D_W$ be a weaving motif of an oriented periodic weave $W_{\infty}$ and let $r$ be a positive integer. The $r$-parallel of $D_W$, denoted $(D_W)^r$, is a weaving motif in which each component of $D_W$ is replaced by $r$ parallel copies that follow the same crossing information as the original component. 
\end{definition}

\begin{lemma}\label{lem:5-6}
Let $D_W$ be a weaving motif of an oriented periodic weave $W_{\infty}$ and $(D_W)^r$, the $r$-parallel of $D_W$.
If $D_W$ is A-adequate (resp. B-adequate), then $(D_W)^r$ is A-adequate (resp. B-adequate).
\end{lemma}

The proof is immediate following \cite{Boden2022} by observing that the state $S'_A$ of $(D_W)^r$ consist of $r$ parallel copies of each curve of the state $S_A$ of $D_W$. 
So if each component of $S_A$ never form a self-crossing at a former crossing of $D_W$, then we have the same conclusion for $S'_A$. \\

\subsection{Relation between the number of crossings and the writhe}\label{sec:5-2}

~

~

In this section, a generalization of Proposition.~\ref{prop:3-8}  gives us the following lemma about the degree of the bracket polynomial.

\begin{lemma}\label{lem:5-7}
Let $\max \mathrm{deg} (\langle D_W \rangle$ and $\min \mathrm{deg}  (\langle D_W \rangle)$ be respectively the maximal and the minimal degree of the bracket polynomial $\langle D_W \rangle$ of a given weaving motif $D_W$. Let $S_A$ denote the state of $D_W$ in which all crossings are $A$-smoothed, and let $S_B$ denote the state of $D_W$ in which all crossings are $B$-smoothed. Let $C$ be the number of crossings in $D_W$. Then,
\begin{equation}
\begin{aligned}
\max \mathrm{deg} (\langle D_W \rangle) &\leq C + 2c_{S_A}, \mbox{ with equality if $D_W$ is A-adequate, }\\
\min \mathrm{deg}  (\langle D_W \rangle) \, &\geq -C - 2c_{S_B}, \mbox{ with equality if $D_W$ is B-adequate. } 
\end{aligned}
\end{equation}
\end{lemma}

Moreover, we can also state the following key lemma which brings together the number of crossings in an A-adequate weaving motif and its writhe, while also connecting the linking numbers to the $r$-parallels. 

\begin{lemma}\label{lem:5-8}
Let $D_{W_1}$ and $D_{W_2}$ be two weaving motifs defined as the quotient of an oriented periodic weave $W_{\infty}$ by a fixed point lattice $\Lambda$, with $C_1$ and $C_2$ crossings, respectively. Suppose that $D_{W_1}$ is A-adequate. Let $w_r(D_{W_1})$ and  $w_r(D_{W_2})$ denote the writhes of $D_{W_1}$ and $D_{W_2}$, respectively. Then, 
$$C_1 - w_r(D_{W_1}) \leq C_2 - w_r(D_{W_2}).$$
\end{lemma}

\begin{proof}
We start by indexing each of the components of $D_W$ such that a thread $t_i$ of $W_{\infty}$ is mapped to the curves $D^i_{W_1}$ and $D^i_{W_2}$ in $D_{W_1}$ and $D_{W_2}$ respectively. Since $D_{W_1}$ is plus-adequate, it does not admit any self-crossing by definition. However, $D_{W_2}$ can contain finitely many self-crossings with sign $\pm 1$. 
Nevertheless, it is always possible to choose an integer $k_i$ that cancels the writhe of the component $D^i_{W_2}$ containing self-crossings. In other words, by performing appropriate Reidemeister moves of type $R_1$, one can add $k_i$ twists of type $+1$ if the original writhe of $D^i_{W_2}$ is negative, or of type $-1$ if it is positive. Therefore, for all integer $i$, we have,
$$w_r(D^i_{W_1}) = w_r(D^i_{W_2}) + k_i = 0. $$ 

By performing these moves to $D^i_{W_2}$, we created $\sum_{i \geq 1} k_i$ self-crossings to $D_{W_2}$. This new diagram is denoted by $D'_{W_2}$. Now to compare $w_r(D'_{W_1})$ with $w_r(D'_{W_2})$, one needs also to consider the crossings which involve distinct components, and thus their linking numbers. 

Therefore, we have,

\begin{equation}
\begin{aligned}
w_r(D_{W_1}) &= \sum_{i \geq j} lk(t_i, t_j)_{\Lambda} \\
w_r(D'_{W_2}) &= \sum_{i \geq1} w_r(D^i_{W_2}) + \sum_{i \geq 1} k_i + \sum_{i \geq j}  lk(t_i, t_j)_{\Lambda} = \sum_{i \geq j}  lk(t_i, t_j)_{\Lambda}
\end{aligned}
\end{equation}

Indeed, since the linking numbers are invariant for a fixed lattice, we obviously have $w_r(D_{W_1}) = w_r(D'_{W_2})$.

We now consider the $r$-parallel $(D_{W_1})^r$ and $(D'_{W_2})^r$. They are both constructed from equivalent weaving motifs for the same point lattice by adding $r - 1$ parallel components. Therefore, they are equivalent by definition and have the same bracket polynomial. Moreover, since every crossing of $D_{W_1}$ and $D'_{W_2}$  corresponds to $r^2$ crossings of $(D_{W_1})^r$ and $(D'_{W_2})^r$, we see that 
$$w_r((D_{W_1})^r) = r^2 w_r(D_{W_1}) = r^2 w_r(D'_{W_2}) = w_r((D'_{W_2})^r).$$

By the definition of the bracket polynomial, it follows that, 
$$\langle (D'_{W_1})^r \rangle = \langle (D'_{W_2})^r \rangle .$$

Let $c_{S^1_A}$ (resp. $c_{S^2_A}$) denote the number of null-homotopic components in the state $S_A$ of $D_{W_1}$ (resp. $D_{W_2}$). Adding self-crossings to $D_{W_2}$ means that the number of connected components in the state $S_A$ of $D'_{W_2}$ becomes $c_{S^2_A} + \sum_{i \geq 1} k_i$. Then, when we pass to the $r$-parallels, we find that the number of connected components in the state $S_A$ of $(D_{W_1})^r$ and $(D'_{W_2})^r$ becomes $r (c_{S^1_A})$ and $r (c_{S^2_A} + \sum_{i \geq 1} k_i)$, respectively.  

Moreover, adding self-crossings in $D_{W_2}$ means that we increase the number of crossings in $D'_{W_2}$, which becomes $C_2 + \sum_{i \geq 1} k_i$. Furthermore, making $r-$parallels means that the number of crossings in $(D_{W_1})^r$ and $(D'_{W_2})^r$ becomes $r^2 C_1$ and $r^2 (C_2 + \sum_{i \geq 1} k_i)$ respectively. Since $D_{W_1}$ is A-adequate, we have, by Lemma.~\ref{lem:5-6} , that $(D_{W_1})^r$ is also A-adequate.  

Thus, from Lemma.~\ref{lem:5-7} , we conclude that, 
$$\max \mathrm{deg} (\langle (D'_{W_1})^r \rangle) = r^2 C_1 + 2r c_{S^1_A} $$                                                                                            
and,
$$\max \mathrm{deg} (\langle (D'_{W_2})^r \rangle) = r^2(C_2 + \sum_{i \geq 1} k_i) + 2r (c_{S^2_A} + \sum_{i \geq 1} k_i)  $$                                                                                           
Since $\langle (D_{W_1})^r \rangle = \langle (D'_{W_2})^r \rangle$,  
then, $\max \mathrm{deg} (\langle (D_{W_1})^r \rangle) = \max \mathrm{deg} (\langle (D'_{W_2})^r \rangle)$, 
and thus, for all positive integers $r,$ we have
$$r \, C_1 + 2 c_{S^1_A} \leq  r (C_2 + \sum_{i \geq 1} k_i)  (c_{S^2_A} + \sum_{i \geq 1} k_i)$$
Therefore,
$$C_1 \leq C_2 + \sum_{i \geq 1} k_i. $$
And since for all positive integer $i$, we have $w_r(D^i_{W_2}) + k_i = 0$, then,
$$C_1 \leq C_2 - \sum_{i \geq 1} w_r(D^i_{W_2}).$$
Again, since the linking number is an invariant, we have
$$lk(D^i_{W_1},D^j_{W_1}) = lk(D^i_{W_2},D^j_{W_2})$$
So as desired, we finally have from (5.2) that,
$$C_1 - w_r(D_{W_1}) \leq C_2 - w_r(D_{W_2}).$$
\end{proof}

\subsection{Tait's Second Conjecture}\label{sec:5-3}

~
~

We end this paper with the statement and proof of the Tait’s Second Conjecture for periodic weaves, as follows.

\begin{theorem}\label{th:3-14} \textbf{(Tait’s Second Conjecture for weaves)}
Any two minimal reduced alternating weaving motifs of an oriented periodic alternating weave have the same writhe. 
\end{theorem}

\begin{proof}
    Let $D_{W_1}$ and $D_{W_2}$ be two minimal reduced alternating weaving motifs of the same oriented periodic weave $W_{\infty}$, which are therefore also adequate. Let $C_1$ and $C_2$ denote the number of crossings in $D_{W_1}$ and $D_{W_2}$, respectively. Then, by the previous Lemma, we have $C_1 - w_r(D_{W_1}) \leq C_2 - w_r(D_{W_2})$,  and $C_2 - w_r(D_{W_2}) \leq C_1 - w_r(D_{W_1})$, and thus, $C_1 - w_r(D_{W_1}) = C_2 - w_r(D_{W_2})$. Moreover, such particular weaving diagrams have the same crossing number from Tait’s First Conjecture, $C_1 = C_2$, and therefore, $w_r(D_{W_1}) = w_r(D_{W_2})$. This finally proves Tait's Second Conjecture for periodic alternating weaves.
\end{proof}


%

\end{document}